\documentclass[11pt]{article}
\usepackage{amsfonts,latexsym,rawfonts,amsmath,amssymb,amsthm,graphicx}
\numberwithin{equation}{section}
\newtheorem{theorem}{Theorem}[section]
\newtheorem{lem}[theorem]{Lemma}
\newtheorem{thm}[theorem]{Theorem}
\newtheorem{pro}[theorem]{Proposition}

\def\s{\,\,\,\,}

\def\endproof{$\hfill\Box$\\}

\title{A remark on a Bernstein type theorem for entire Willmore graphs in $R^3$}
\author{Yong Luo\footnote{The author is supported by the DFG Collaborative Research Center SFB/TR71.}, Jun Sun\footnote{The author was supported by NSF in China, No. 11001268, No. 11071236.}}
\date{}
\begin{document}
\maketitle
\begin{abstract}
In this note we prove that every two-dimensional entire Willmore graph in $R^3$ with square integrable mean curvature is a plane.

\end{abstract}

\section{Introduction}
For a two dimensional closed surface $\Sigma$ and an immersion $f:\Sigma\rightarrow R^3$ the Willmore functional is defined by
\begin{eqnarray}
W(f)=\frac{1}{4}\int_\Sigma H^2d\mu_g,
\end{eqnarray}
where $g$ is the induced metric, $H$ is the mean curvature of $\Sigma$ and $d\mu_g$ is the area element.

Critical points of $W$ are called Willmore surfaces and they are solutions of the Euler-Lagrangian equation
\begin{eqnarray}\label{E-L equation}
\Delta_gH+\frac{1}{2}H^3-2HK=0,
\end{eqnarray}
where $\Delta_g$ is the Laplace-Beltrami operator of the induced metric and $K$ is the Gauss curvature of $\Sigma$. Moreover we denote by $A$ the second fundamental form of the immersion $f$.

It is obviously that every minimal surfaces in $R^3$ is a Willmore surface. The classic Bernstein theorem says that every entire minimal graph in $R^3$ is a plane. It is nature to ask whether there are such type results for entire Willmore graphs. Recently, Chen-Lamm (\cite{CL}) started the study of this direction and proved that:
\begin{thm}{(Chen-Lamm)}\label{Chen-Lamm}
Every smooth, entire graphical solution of (\ref{E-L equation}) with finite $L^2$ norm of the second fundamental form is a plane.
\end{thm}
 In the same paper, Chen-Lamm (\cite{CL}) proposed a question that whether we can weaken the assumption $A\in L^2$ to be $H\in L^2$. In this note, we will answer this question in \cite{CL}:
\begin{theorem}\label{Luo-Sun}
Every smooth, entire graphical solution of (\ref{E-L equation}) with finite $L^2$ norm of the mean curvature is a plane.
\end{theorem}
This is also a generalization of the classic Bernstein theorem  since minimal entire graphs obviously satisfy the assumptions of the above theorem.

Actually we will prove that the total Gauss curvature of a smooth entire graph in $R^3$ with square integrable mean curvature is zero:
\begin{pro}\label{total curvature}
Let $\Sigma=(x,y,u(x,y))$, where $u:R^2\rightarrow R$, be a smooth graph with square integrable mean curvature. Then
$$\int_\Sigma Kd\mu_g=0.$$
\end{pro}
Note that $|A|^2=H^2-2K$. Hence for an entire graph in $R^3$ with square integral mean curvature, we have
\begin{eqnarray}
\int_\Sigma|A|^2d\mu_g=\int_\Sigma H^2d\mu_g.
\end{eqnarray}
Then Theorem \ref{Luo-Sun} follows directly from Theorem \ref{Chen-Lamm}.

\textbf{Acknowledgment} The first author would like to thank his advisor, Professor Guofang Wang, for bringing his attention to the paper \cite{CL} and discussions.

\section{Proof of Proposition \ref{total curvature}}
\textbf{Proof of Proposition \ref{total curvature}:}
For the convenience we assume that $0\in\Sigma=(x,y,u(x,y))$. At first we will show that every entire graph in $R^3$ with square integral mean curvature has quadratic area growth by following the ideas of \cite{CL}:
\begin{lem}
Let $0\in\Sigma=(x,y,u(x,y))$, where $u:R^2\rightarrow R$, be a smooth entire graph with square integrable mean curvature. Then $\Sigma$ has quadratic area growth.
\end{lem}
\proof Recall that in the proof of Proposition 2.1 of \cite{CL}, using a calibration argument (see \cite{CM}) Chen and Lamm proved the following inequality which holds for almost every $R>0$:
\begin{eqnarray}\label{e2.1}
|\Sigma\cap B(R)|\leq2\pi R^2+\int_{K(R)}|H(x,y,u(x,y))|dxdydz,
\end{eqnarray}
where $K(R)$ is a 3-dimensional subset of $B(R)$.

In order to show the quadratic area growth of $\Sigma$, we only need to control the second term of the right hand side. It is easy to see that (see page 1 of \cite{CM}) the induced area element on $\Sigma$ is given by $d\mu_g=\sqrt{1+|Du|^2}$, where $Du=(u_x,u_y)$. Then by H\"older inequality,
\begin{eqnarray}\label{e2.3}
   & & \int_{K(R)}|H(x,y,u(x,y))|dxdydz \nonumber\\
   &\leq & \int_{B(R)}|H(x,y,u(x,y))|dxdydz
   \leq  \int_{-R}^{R}\int_{x^2+y^2\leq R^2}|H(x,y,u(x,y))|dxdydz \nonumber\\
   &=& 2R\int_{x^2+y^2\leq R^2}|H(x,y,u(x,y))|dxdy
   \leq 2R\left(\int_{x^2+y^2\leq R^2}|H(x,y,u(x,y))|^2dxdy\right)^{\frac{1}{2}}(\pi R^2)^{\frac{1}{2}}\nonumber\\
   &\leq& 2\sqrt{\pi}R^2\left(\int_{x^2+y^2\leq R^2}|H(x,y,u(x,y))|^2\sqrt{1+|Du|^2}dxdy\right)^{\frac{1}{2}}\nonumber\\
   &\leq& 2\sqrt{\pi}R^2\left(\int_{\Sigma}|H|^2d\mu_g\right)^{\frac{1}{2}}\leq C_1R^2,
\end{eqnarray}
where the last inequality used our assumption.
By (\ref{e2.1}) and (\ref{e2.3}), we know that
\begin{equation}\label{e2.4}
    |\Sigma \cap B(R)|\leq (2\pi+C_1) R^2\equiv C_2R^2.
\end{equation}
Since both sides of this estimate are continuous in $R$, we actually get this estimate for all $R>0$. This proves that $\Sigma$ has quadratic area growth.
\endproof

Next we prove that the total Gauss curvature vanishes. We fix an orientation for $\Sigma$ and we denote by $n:\Sigma\rightarrow S^2$ the Gauss map of the surface $\Sigma$. Furthermore let $\beta$ be the area 2-form on the unit sphere $S^2$. As $\Sigma$ is a graph, $n(\Sigma)$ lies in the upper hemisphere $S_+^2$ which is contractible, thus $\beta=d\alpha$ for some smooth 1-form $\alpha$ on $S_+^2$. The Gauss curvature $K$ of $\Sigma$ measures the area distortion under the Gauss map:
$$Kd\mu_g=K\sqrt{g}dx\wedge dy=\det(dn)dx\wedge dy=n^\ast\beta=n^\ast d\alpha=d(n^\ast\alpha).$$
For the smooth 1-form $\alpha$ there is a constant $C_\alpha$ depending on $\alpha$ such that
$$|n^\ast\alpha|\leq C_\alpha|dn|=C_\alpha|A|.$$
Then let $0\leq\eta\leq1$ be a smooth function with compact support and we calculate
\begin{eqnarray}\label{ES}
|\int_\Sigma\eta^2K\sqrt{g}dx\wedge dy|=|\int_\Sigma2\eta d\eta\wedge n^\ast\alpha|\leq4C_\alpha(\int_\Sigma\eta^2|A|^2d\mu_g)^\frac{1}{2}(\int_\Sigma|\nabla\eta|^2d\mu_g)^\frac{1}{2}.
\end{eqnarray}
Now by the Gauss equation, we have$|A|^2=H^2-2K$, thus
\begin{eqnarray*}
\int_\Sigma\eta^2|A|^2d\mu_g&=&\int_\Sigma\eta^2 H^2d\mu_g-2\int_\Sigma \eta^2Kd\mu_g
\\&\leq& C_3+2|\int_\Sigma \eta^2Kd\mu_g|
\\&\leq&8C_\alpha(\int_\Sigma\eta^2|A|^2d\mu_g)^\frac{1}{2}(\int_\Sigma|\nabla\eta|^2d\mu_g)^\frac{1}{2}+C_3.
\end{eqnarray*}
Solving this inequality yields
$$\int_\Sigma\eta^2|A|^2d\mu_g\leq8C_\alpha(\int_\Sigma |\nabla\eta|^2d\mu_g)^\frac{1}{2}+\sqrt{C_3}.$$
Putting the above inequality to (\ref{ES}) we get
$$|\int_\Sigma\eta^2K\sqrt{g}dx\wedge dy|\leq4C_\alpha(8C_\alpha(\int_\Sigma|\nabla\eta|^2d\mu_g)^\frac{1}{2}+\sqrt{C_3})(\int_\Sigma|\nabla\eta|^2d\mu_g)^\frac{1}{2}.$$

Now it is standard to prove that the total Gauss curvature vanishes. For the convenience of the reader we provide the proof in the sequel.

For every $\sigma>1$ we define a radially symmetric function $\eta_\sigma$ on $R^3$ as follows:
$$\eta_\sigma=1\s for\s|x|\leq\sqrt{\sigma};\s \eta_\sigma=2-2\frac{\log|x|}{\log\sigma}\s for\s\sqrt{\sigma}\leq|x|\leq\sigma;\s\eta_\sigma=0\s for\s|x|>\sigma.$$
For $|\nabla|x||\leq|D|x||=1$, we have
$$|\nabla\eta_\sigma|\leq\frac{2}{|x|\log\sigma}.$$
Hence for $\sigma$ large enough, we have
\begin{eqnarray*}
\int_\Sigma|\nabla\eta_\sigma d\mu_g|\leq\frac{4}{\log^2\sigma}\sum_{\frac{1}{2}\log\sigma}^{\log\sigma-1}\int_{\Sigma\cap B(e^{l+1})\setminus\Sigma\cap B(e^l)}|x|^{-2}d\mu_g
\leq\frac{4}{\log^2\sigma}\sum_{\frac{1}{2}\log\sigma}^{\log\sigma-1}C_4e^2
\leq\frac{C_5}{\log\sigma},
\end{eqnarray*}
where in the second inequality we have used (\ref{e2.4}).

Let $\sigma\to\infty$ we get
\begin{eqnarray*}
|\int_\Sigma Kd\mu_g|&=&\lim_{\sigma\to\infty}|\int_\Sigma\eta_\sigma^2Kd\mu_g|
\leq\lim_{\sigma\to0}4C_\alpha(8C_\alpha(\int_\Sigma|\nabla\eta_\sigma|^2d\mu_g)^\frac{1}{2}+\sqrt{C_3})(\int_\Sigma|\nabla\eta_\sigma|^2d\mu_g)^\frac{1}{2}
=0,
\end{eqnarray*}
which completes the proof of proposition \ref{total curvature}.

\section{Appendix}
Let $\Sigma=(x,y,u(x,y))$ be a Willmore graph in $R^3$ where $u:R^2\rightarrow R$. Standard calculations then yield that
\begin{eqnarray*}
H&=&div(\frac{\nabla u}{v}), K=\frac{\det\nabla^2u}{v^4} and
\\ \Delta_g H&=&\frac{1}{v}div((vI-\frac{\nabla u\otimes\nabla u}{v})\nabla H)
\end{eqnarray*}
where $v=\sqrt{1+|Du|^2}$. From the calculations in \cite{DD} we then get that the Willmore equation (\ref{E-L equation}) can be rewritten as follows:
\begin{eqnarray}
div(\frac{1}{v}((I-\frac{\nabla u\otimes\nabla u}{v^2})\nabla(vH)-\frac{1}{2}H^2\nabla u)).
\end{eqnarray}

 \textbf{Notes:} In the paper of \cite{CL}, the graph condition is used to derive the fact that the image of Gauss map is contained in the upper hemi-sphere, which is contractible, and to exclude the case that the density of the entire Willmore graph at infinity is 2. If the graph function depends only on one variable, they(\cite{CL}, Appendix) proved the theorem without the finiteness condition on the integral of the square of the second fundamental form, by a simple calculation toward the above graphical equation, but in the general case, the graphical equation is not used. Hence we are nature to guess that there should be another different proof of this theorem by using the graphical equation essentially, possibly with weaker conditions.

{}

\vspace{1cm}\sc

Yong Luo

Mathematisches Institut,

Albert-Ludwigs-Universit\"at Freiburg,

Eckerstr. 1, 79104 Freiburg, Germany.

{\tt yong.luo@math.uni-freiburg.de}\\

Jun Sun

Math. Group,

The Abdus Salam International Centre for Theoretical Physics

Trieste 34100, Italy.

{\tt jsun@@ictp.it}
\end{document}